\documentclass{amsart}

\usepackage{amsmath}
\usepackage{amscd}
\usepackage{amssymb}

\newcommand{\cal}{\mathcal}

\newcommand{\bC}{{\Bbb C}}

\newcommand{\bZ}{{\Bbb Z}}

\newcommand{\cT}{{\cal T}}

\DeclareMathOperator{\ch}{ch}

\DeclareMathOperator{\Img}{Im}
\DeclareMathOperator{\Ker}{Ker}

\newtheorem{theorem}{Theorem}[section]
\newtheorem{theorem/definition}{Theorem/Definition}[section]

\newtheorem{proposition}{Proposition}[section]
\newtheorem{lemma}{Lemma}[section]

\newtheorem{corollary}{Corollary}[section]

\theoremstyle{remark}
\newtheorem{remark}{Remark}[section]

\theoremstyle{definition}
 \newtheorem{example}{Example}[section]

\begin{document}
\title[Calculations of Hirzebruch $\chi_y$ genera of symmetric products]
{Calculations of the Hirzebruch $\chi_y$ genera of 
symmetric products by the holomorphic Lefschetz formula}
\author{Jian Zhou}
\address{Department of Mathematics\\
Texas A\&M University\\
College Station, TX 77843}
\email{zhou@math.tamu.edu}
\begin{abstract}
We calculate the Hirzebruch $\chi_y$ and $\hat{\chi}_y$-genera 
of symmetric products of closed complex manifolds 
by the holomorphic Lefschetz formula  of Atiyah and Singer \cite{Ati-Sin}.
Such calculation rederive some formulas proved in an earlier paper \cite{Zho}
by a different method.
\end{abstract}

\maketitle

Let $M$ be a smooth manifold, and $G$ a finite group of diffeomorphisms.
There are two kinds of interesting cohomology theories 
for the orbifold $M/G$.
The first is the orbifold de Rham cohomology $H^*(M/G)$ 
introduced in \cite{Sat},
for which we easily see to have an isomorphism
\begin{eqnarray} \label{orbifoldcohomology}
H^*(M/G) \cong H^*(M)^G.
\end{eqnarray}
The second is the delocalized equivariant cohomology 
\begin{eqnarray} \label{Delocalized} \;\;\;\;\;\;
H^*(M, G) = \left(\bigoplus_{g \in G} H^*(M^g)\right)^G
= \bigoplus_{[g] \in G_*} H^*(M^g)^{Z_g}
= \bigoplus_{[g] \in G_*} H^*(M^g/Z_g),
\end{eqnarray}
introduced by 
Baum-Connes \cite{Bau-Con} in the study of the equivariant $K$-theory.
Here $G_*$ denotes the set of conjugacy classes of $G$
and $Z_g$ denotes the centralizer of $g$.
For each of these cohomology theories,
one can define an Euler number,
denoted by $\chi(M/G)$ and $\chi(M, G)$ respectively.
The latter first appeared in the string theory on orbifolds
in a different form (cf. Dixon-Harvey-Vafa-Witten \cite{Dix-Har-Vaf-Wit}):
\begin{eqnarray} \label{def:orbifoldEuler}
\chi(M, G) 
= \frac{1}{|G|} \sideset{}{'}\sum_{g, h} \chi(M^{\langle g, h \rangle}),
\end{eqnarray}
where $\langle g, h \rangle$ is the group generated by $g$ and $h$,
the sum is taken over commutating pairs $(g, h) \in G \times G$.
From Atiyah-Segal \cite{Ati-Seg} and Hirzebruch-H\"{o}fer \cite{Hir-Hof}
one knows that $\chi(M, G)$ as defined in (\ref{def:orbifoldEuler})
is the Euler number of $H^*(M, G)$ or equivalently $K_G^*(M)$,
since one can easily show that
\begin{eqnarray} \label{reducetofix2}
\chi(M, G) = \sum_{[g] \in G_*} \chi(M^g/Z_g).
\end{eqnarray}
A classical analogue of this formula is of course the Lefschetz formula:
\begin{eqnarray} \label{reducetofix1}
\chi(M/G) = \frac{1}{|G|}\sum_{g \in G} \chi(M^g) 
= \sum_{[g] \in G_*} \frac{1}{|Z_g|} \chi(M^g).
\end{eqnarray}

A very interesting class of examples 
are provided by the symmetric products $X^{(n)} = X^n/S_n$ of a manifold $X$.
The following two formulas have been proved by various authors:
\begin{eqnarray} 
&& \sum_{n \geq 0} \chi(X^{(n)}) q^n = \frac{1}{(1 - q)^{\chi(X)}}, 
	\label{gfEuler1}\\
&&  \sum_{n \geq 0} \chi(X^n, S_n) q^n 
= \prod_{l \geq 1} \frac{1}{(1 - q^l)^{\chi(X)}}. \label{gfEuler2}
\end{eqnarray}
There are two approaches to such formulas.
Macdonald \cite{Mac} proved (\ref{gfEuler1}) using the isomorphism
\begin{eqnarray} \label{eqn:isomorphism1}
H^*(X^{(n)}) \cong H^*(X^n)^{S_n} \cong S^n(H^*(X)),
\end{eqnarray}
and Vafa and Witten \cite{Vaf-Wit1} proved (\ref{gfEuler2})  
using the isomorphism
\begin{eqnarray} \label{eqn:isomorphism2}
\oplus_{n \geq 0} H^*(X^n, S_n) \cong S(\oplus_{n \geq 0} H^*(X)).
\end{eqnarray}
On the other hand,
in Zagier \cite{Zag}, \S9,
one can find a proof of (\ref{gfEuler1}) by (\ref{reducetofix1}),
while in Hirzebruch-H\"{o}fer \cite{Hir-Hof} a proof of (\ref{gfEuler2})
by (\ref{reducetofix2}).

It is interesting to consider the complex version of (\ref{gfEuler1})
and (\ref{gfEuler2}) by both approaches.
The first approach has been carried in our earlier paper \cite{Zho},
we will carry out the second approach in this paper.
For a closed complex manifold $M$,
a complex analogue of the Euler number is the Hirzebruch $\chi_y$-genus:
\begin{eqnarray*}
\chi_y(M) = \sum_{p, q \geq 0} (-1)^q h^{p, q}(M) y^p.
\end{eqnarray*}
For a closed complex orbifold $M/G$,
define $\chi_y(M/G)$ and $\chi_y(M, G)$ using the Hodge numbers of
two versions of Dolbeault cohomology theories 
$$H^{*, *}(M/G) = H^{*, *}(M)^G$$ 
and
$$H^{*, *}(M, G) = \bigoplus_{[g] \in G_*} H^{*, *}(M^g)^{Z_g} =
\bigoplus_{[g] \in G_*} H^{*, *}(M^g/Z_g)$$
respectively.
For the latter,
as in \cite{Zho},
we use a graded shift suggested by physicists:
to each connected component of $M^g$,
if the eigenvalues of $g$ on the normal bundle are
$\exp (\sqrt{-1}\theta_1), \cdots, \exp (\sqrt{-1}\theta_r)$,
set 
$$F_g = \frac{1}{2\pi}(\theta_1 + \cdots + \theta_r).$$
The Hodge numbers of $H^{*, *}(M, G)$ are given by
\begin{eqnarray} \label{Hodgenumber}
h^{p, q}(M, G) = \sum_{[g] \in G_*} h^{p-F_g, q-F_g}(M^g/Z_g).
\end{eqnarray}
For a closed complex manifold $X$,
the complex versions of (\ref{eqn:isomorphism1}) and (\ref{eqn:isomorphism2}) 
are the following isomorphisms respectively:
\begin{eqnarray} 
&& H^{*, *}(X^{(n)}) \cong S^n(H^{*, *}(X)), \label{eqn:complexisomorphism1} \\
&& \bigoplus_{n \geq 0} H^{*, *}(X^n, S_n) \cong 
S(\bigoplus_{n \geq 0} H^{*, *}(X)[F_g, F_g]) \label{eqn:complexisomorphism2}
\end{eqnarray}
respectively.
They have been used in \cite{Zho} to prove the following 
analogues of (\ref{gfEuler1}) and (\ref{gfEuler2}) respectively:
\begin{eqnarray} 
&& \sum_{n \geq 0} \chi_{-y}(X^{(n)}, S_n) q^n 
= \exp \left( \sum_{l \geq 1} \frac{\chi_{-y^l}(X)}{l} q^l \right), 
\label{gfchiy1} \\
&& \sum_{n \geq 0} \chi_{-y}(X^n, S_n) q^n 
= \exp \left( \sum_{l \geq 1} 
\frac{q^l}{l}
\frac{\chi_{-y^l}(X)}{1 - (y^{\dim X/2}q)^l}  \right). \label{gfchiy2}
\end{eqnarray}
In this paper, 
we use the complex analogue of (\ref{reducetofix1}),
the holomorphic Lefschetz formula of Atiyah-Singer \cite{Ati-Sin},
to rederive (\ref{gfchiy1}).
The following complex analogue of (\ref{reducetofix2}) 
follows easily from (\ref{gfchiy1}):
\begin{eqnarray} \label{complexreducetofix2}
\chi_y(M, G) = \sum_{[g] \in G_*} (-y)^{F_g} \chi_y(M^g/Z_g).
\end{eqnarray}
We will use it to derive (\ref{gfchiy2}).
Incidentally,
we find a description of the Adams operation in terms of 
the localized theorem proved by Atiyah and Segal \cite{Ati-Seg}.
See Proposition \ref{prop:Adams}.
As in \cite{Zho},
we also have a version for $H^{-*, *}$ and $\hat{\chi}_y$.

\section{Preliminaries}

\subsection{Equivariant $K$-theory}
\label{sec:equivariantK}

Recall that for a $G$-manifold $M$,
$K^0_G(M)$ is the abelian group generated by the complex $G$-vector bundles,
while $K_G^1(M)$ can defined as the kernel of the restriction map
$$K_G^0(M \times S^1) \to K_G^0(M)$$
given by the inclusion of point in $S^1$.
Atiyah and Segal \cite{Ati-Seg} proved the following the following

\begin{theorem}
There is a natural isomorphism
$$\phi = \bigoplus_{[g] \in G_*} \phi_g: K_G^*(M) \otimes \bC 
\cong \bigoplus_{[g] \in G_*} [K^*(M^g) \otimes \bC]^{Z_g}.$$
\end{theorem}

The isomorphism on $K^0_G(M)$ can be explicitly given as follows.
If $E$ is a $G$-vector bundle over $M$,
its restriction to $M^g$ is acted on fiberwise by $g$ and 
so decompose as a direct sum of subbundles $E_{\xi}$ 
for each eigenvalue $\xi$ of $g$.
Then 
$$\phi_g(E) = \sum_{\xi} \xi E_{\xi}.$$ 

\subsection{Tensor products and Adams operations}

Given a vector bundle $\pi: E \to X$,
$E^{\otimes n}$ is a vector bundle on $X$.
There is a natural $S_n$-action on $E^{\otimes n}$ given as follows:
$$\sigma(v_1 \otimes \cdots \otimes v_n) 
= v_{\sigma^{-1}(1)} \otimes \cdots \otimes v_{\sigma^{-1}(n)},$$
where $\sigma \in S_n$, $v_1, \cdots, v_n \in E_x = \pi^{-1}(x)$, $x \in X$.
Atiyah \cite{Ati1} (Proposition 2.2) showed that $E \mapsto E^{\otimes n}$
defines a map 
$$\otimes n: K(X) \to K_{S_n}(X) \cong K(X) \otimes R(S_n).$$
The {\em Adams operations} $\psi^n: K(X) \rightarrow K(X)$ 
is essentially defined as follows:
if $E = L_1 \oplus \cdots \oplus L_r$,
where 
$L_1, \cdots, L_r$ are line bundles,
then
$$\psi^n(E) = L_1^{\otimes n} \oplus \cdots \oplus L_r^{\otimes n}.$$

\begin{proposition} \label{prop:Adams}
For a vector bundle $E$ over $X$,
$\phi_{\sigma_n}(E^{\otimes  n}) = \psi^n(E)$.
\end{proposition}

\begin{proof}
By splitting principle,
we may assume that $E = \oplus_{m = 1}^r L_m$,
where $L_m$'s are line bundles.
Then 
$$E^{\otimes n} = \bigoplus_{1 \leq j_1, \cdots, j_n \leq r}  
L_{j_1} \otimes \cdots \otimes L_{j_n}.$$
Consider the index set
$\{ (j_1, \cdots, j_n): \; 1 \leq j_1, \cdots, j_n \leq r \}$.
The action of $\sigma_n$ or equivalently the cyclic group $\bZ_n$
on $E^{\otimes n}$ corresponds to the cycling of the indices:
$$\sigma_n (j_1, \cdots, j_n) = (j_n, j_1, \cdots, j_{n-1}).$$ 
For each $J = (j_1, \cdots, j_n)$,
let $L_J = L_{j_1} \otimes \cdots \otimes L_{j_n}$,
and $V_J$ be the subbundle spanned by
line bundles $L_{\sigma_n^k(J)}$, $k=0, \cdots, n-1$.
Then $V_J$ is invariant under the action of $\bZ_n = \langle \sigma_n \rangle$.
Since $E^{\otimes n}$ is a direct sum of such $V_J$'s, 
it suffices to find $\phi_{\sigma_n}(V_J)$ for all $J$.
When $j_1 = \cdots = j_n = m$,
$V_J = L_m^{\otimes n}$ is fixed by $\sigma_n$.
In other words,
$\sigma_n$ has eigenvalue $1$ on $\oplus_{m=1}^r L_m^{\otimes n} = \psi^n(E)$.
When $j_1, \cdots, j_n$ are not all identical,
there are two cases to consider.
If the orbit of $J$ has length $n$,
then $V_J$ is the direct sum of $n$-copies of $L_J$,
and each fiber is a regular representation of 
$\bZ_n = \langle \sigma_n \rangle$.
Therefore,
$\sigma_n$ has eigenvalues $e^{2k\pi\sqrt{-1}}/n$, 
$k = 0, \cdots, n-1$.
Hence 
$$\phi_{\sigma_n}(V_J) = \sum_{k = 0}^{n-1} e^{2k \pi \sqrt{-1}/n} L_J = 0.$$
Another case is that the orbit of $J$ has length $1 < l < n$.
Then it is easy to see that $l | n$, $j_k = j_{k + l}$, 
i.e., 
$J$ is of the form 
$$J = (j_1, \cdots, j_l, j_1, \cdots, j_l, \cdots, j_1, \cdots, j_l).$$
$V_J$ is isomorphic to the direct sum of $l$ copies of 
$(L_{j_1} \otimes \cdots \otimes L_{j_l})^{n/l}$.
Furthermore,
$\bZ_n$ acts on $V_J$ via $\bZ_l = \bZ_n /\bZ_{n/l}$-action 
which correpsonds to $(j_1, \cdots, j_l) \mapsto (j_n, j_1, \cdots, j_{l-1})$.
Therefore,
one sees that $\sigma_n$ has eigenvalues $e^{2k\pi\sqrt{-1}/l}$, 
$k = 0, \cdots, n-1$.
Hence 
$$\phi_{\sigma_n}(V_J) = \sum_{k=0}^{k-1} e^{2k\pi\sqrt{-1}/l}
(L_{j_1} \otimes \cdots \otimes L_{j_l})^{n/l}=0.$$
This completes the proof.
\end{proof}

\begin{remark} \label{remark:character}
In the calculation for $j_1, \cdots, j_n$ not all identical,
one can also compute the character of the $\bZ_n$ or $\bZ_l$ representation 
on the fibers of $V_J$.
This again gives the eigenvalues of $\sigma_n$ on $V_J$.
This argument can be generalized to the graded version
(cf. Proposition \ref{prop:gradedAdams}).
\end{remark}

\subsection{External tensor products}

For $i =1, \cdots, n$,
let $E_i \to X_i$ be vector bundles over manifolds $X_i$.
$p_i: X_1 \times \cdots \times X_n \to X_i$ the projections 
onto the $i$-th factors.
The {\em exterior tensor product} of $E_1, \cdots, E_n$ is defined to be:
$$E_1 \boxtimes \cdots \boxtimes E_n 
= p_1^*E_1 \otimes \cdots \otimes p_n^*E_n.$$
For a vector bundle $\pi: E \to X$,
let $E^{\boxtimes n}$ be the external tensor product of $n$-copies of $E$.
There is an $S_n$-action on $X^n$ given by
$$\sigma (x_1, \cdots, x_n) 
= (x_{\sigma^{-1}(1)}, \cdots, x_{\sigma^{-1}(n)}),$$
where $\sigma \in S_n$, $x_1, \cdots, x_n \in X$.
In other words,
each $x_j$ is moved to the $\sigma(j)$-th position.
Similarly, there is an $S_n$-action on $E^{\boxtimes n}$ given by
$$\sigma(v_1 \boxtimes \cdots, v_n) 
= v_{\sigma^{-1}(1)} \boxtimes \cdots \boxtimes v_{\sigma^{-1}(n)},$$
or equivalently,
$$\sigma(p_1^*v_1 \otimes \cdots \otimes p_n^*v_n
= p_1^*v_{\sigma^{-1}(1)} \otimes \cdots \otimes p_n^*v_{\sigma^{-1}(n)},$$
where $\sigma \in S_n$,
$v_j \in E_{x_j} = \pi^{-1}(x_j)$,
$x_j \in X$,
$j = 1, \cdots, n$.
Therefore $E^{\boxtimes n}$ is an $S_n$-bundle over the $S_n$-manifold $X^n$.
Clearly the external tensor product defines a ring homomorphism
(cf. Atiyah \cite{Ati1}, Proposition 3.2):
$$\boxtimes n: K^*(X) \to K_{S_n}(X^n).$$

Let $\Delta_n(X) = \{(x, \cdots, x) \in X^n\}$.
The map $i: X \to \Delta_n(X) \hookrightarrow X^n$ 
given by $x \mapsto (x, \cdots, x)$ induces a homomorphism
$i^*: K_{S_n}(X^n) \to K_{S_n}(X)$.
Clearly $i^*(E^{\boxtimes n}) = E^{\otimes n}$.
Notice that $(X^n)^{\sigma_n} = \Delta_n(X)$. 
From Proposition \ref{prop:Adams},
we get

\begin{corollary} \label{cor:Adams}
For a vector bundle $E$ on $X$,
we have
$i^* \phi_{\sigma_n}(E^{\boxtimes n}) = \psi^n(E)$.
\end{corollary}

\subsection{Riemann-Roch numbers}
Let $\pi: E \to M$ be a holomorphic vector bundle 
over a closed complex $n$-manifold $M$.
Consider the Dolbeault complex 
$$0 \to \Omega^{0, 0}(E) \stackrel{\bar{\partial}_E}{\to} \Omega^{0, 1}(E) 
\stackrel{\bar{\partial}_E}{\to} \cdots \stackrel{\bar{\partial}_E}{\to} 
\Omega^{0, n}(E) \to 0,$$
and the Dolbeault cohomology 
$H^*(M, E) = \Ker \bar{\partial}_E /\Img \bar{\partial}_E$.
The Riemann-Roch number is by definition
$$\chi(M, E) = \sum_{q=0}^n \dim H^q(M, E).$$
The famous Hirzebruch-Riemann-Roch theorem states
\begin{eqnarray*}
\chi(M, E) = \int_M ch(E) \cT (M),
\end{eqnarray*}
where $\cT (M)$ is the Todd class of $M$.

\subsection{Holomorphic Lefschetz theorem} \label{sec:Lefschetz}

For a vector bundle $E$ on $M$,
set 
\begin{eqnarray*}
&& S_t(E) = 1 + t E + t^2 S^2(E) + \cdots, \\
&& \Lambda_t(E) = 1 + t E + t^2 \Lambda^2(E) + \cdots.
\end{eqnarray*}
They are exponential in the sense that
\begin{eqnarray*}
&& S_t(E_1 + E_2) = S_t(E_1) S_t(E_2), \\
&& \Lambda_t(E_1 + E_2) = \Lambda_t(E_1) \Lambda_t(E_2).
\end{eqnarray*}

Now assume that $M$ is a complex $G$-manifold, 
$E$ is a holomorphic $G$-vector bundle.
Denote by $N^g$ the holomorphic normal bundle of a component of $M^g$ in $M$.
Then there is a natural decomposition
$$N^g = \bigoplus_{0 < \theta < 2 \pi} N^g(\theta),$$
where each $N^g(\theta)$ is a holomorphic subbundle on 
which $g$ acts as multiplication by $e^{\sqrt{-1} \theta}$.
Set
$$\chi_g(M, E) = \sum_{q = 0}^n (-1)^q tr (g|_{H^q(M, E)}).$$
The holomorphic Lefschetz theorem states:
\begin{eqnarray} \label{Lefschetz1}
\chi_g(M, E)
= \int_{M^g} \frac{\ch_g (E|_{M^g}) \cT(M^g)}{\ch_g \Lambda_{-1}((N^g)^*)}.
\end{eqnarray}
Recall that if $V$ is a complex 
$G$-bundle on $X$ (with trivial $G$-action on $X$),
one can compute the equivariant characteristic classes as follows:
for any $g \in G$,
$$V = \bigoplus_{0 \leq \theta < 2 \pi} V(\theta),$$ 
where $g$ acts on $V(\theta)$ as multiplication by $e^{\sqrt{-1}\theta}$.
If $\{x_j\}$ denote the Chern roots of $V(\theta)$,
then 
\begin{eqnarray*}
&& \ch_g(V(\theta)) = \sum_j e^{x_j + \sqrt{-1}\theta} 
= e^{\sqrt{-1} \theta} \sum_j e^{x_j} = e^{\sqrt{-1} \theta} \ch V(\theta)
= \ch (e^{\sqrt{-1} \theta} V(\theta)), \\
&& \ch_g \Lambda_{-1} (V(\theta)) = \prod_j (1 - e^{x_j + \sqrt{-1}\theta}), 
\end{eqnarray*}
With the notation of \S \ref{sec:equivariantK},
(\ref{Lefschetz1}) can be rewritten as
\begin{eqnarray} \label{Lefschetz2}
\chi_g(M, E)
= \int_{M^g} \frac{\ch \phi_g (E) \cT(M^g)}{ch_g \Lambda_{-1}((N^g)^*)}.
\end{eqnarray}

\subsection{Index theory on orbifolds} \label{sec:orbifolds}

On the orbifold $M/G$,
$E/G$ is a $V$-vector bundle.
Then one can define orbifold version of Dolbeault operator
$$\bar{\partial}_{E/G}: \Omega^{0, *}(E/G) \to \Omega^{0, *+1}(E/G),$$
where $\Omega^{0, *}(E/G) = \Omega^{0, *}(E)^G$.
One can consider the Dolbeault cohomology 
$$H^*(M/G, E/G) = \Ker \bar{\partial}_{E/G}/\Img \bar{\partial}_{E/G}$$
and the Riemann-Roch number
$$\chi(M/G, E/G) = \sum_{q=0}^n (-1)^q \dim H^q(M/G, E/G).$$
Standard character theory shows that
\begin{eqnarray*}
\chi(M/G, E/G) 
= \frac{1}{|G|} \sum_{g \in G} \chi_g(M, E)
= \sum_{[g] \in G_*} \frac{1}{|Z_g|} \chi_g(M, E).
\end{eqnarray*}
By (\ref{Lefschetz2}),
we then have
$$\chi(M/G, E/G) 
= \sum_{[g] \in G_*} \frac{1}{|Z_g|} 
\int_{M^g} \frac{\ch (\phi_g(E)) \cT_g(M^g)}{\ch_g\Lambda_{-1}((N^g)^*)}$$
String theory on orbifold suggests the consideration of another version:
$$\chi(M, E|G) = \sum_{g \in G_*} \chi(M^g/Z_g, \phi_g(E)/Z_g).$$

\section{Riemann-Roch numbers of symmetric products of vector bundles}

\subsection{Structures of the fixed point sets}
We first recall the structure of $(X^n)^g$ for $g \in S_n$.
Any element of $S_n$ can be uniquely written as a products of
mutually exclusive cycles.
Denote by $N_l(g)$ the number of $l$-cycles in $g$.
The sequence $N(g) = (N_1(g), N_2(g), \cdots)$ is called 
the {\em cycle type} of $g$.
Each cycle type corresponds to a unique conjugacy class,
since permutations with the same cycle type are conjugate to each other.
Given any element $g \in S_n$ of cycle type $N = (N_1, N_2, \cdots)$,
there is an isomorphism 
\begin{eqnarray*}
Z_g \cong S_{N_1} \times (S_{N_2} \wr \bZ^{N_2}_2) \times \cdots 
\times (S_{N_n} \wr \bZ^{N_n}_n),
\end{eqnarray*}
where $S_{N_l}$ corresponds to permutations of the $l$-cycles of $g$,
and 
each $\bZ_l$ corresponds to the cyclic group generated by an $l$-cyle in $g$.
It is clear that
\begin{eqnarray} \label{fix1}
(X^n)^g = \prod_{l=1}^n \Delta_l(X)^{N_l} \cong \prod_{l=1}^n X^{N_l},
\end{eqnarray}
where each $l$-cycle of $g$ contributes a copy of $\Delta_l(X)$.
Since each $Z_l$-factor in $Z_g$ is generated by the corresponding $l$-cycle,
it acts trivially on the corresponding $\Delta_l(X)$.
On the oter hand, 
the $S_{N_l}$-factor acts by permuting the $N_l$ copies of $\Delta_l(X)$.
Therefore
\begin{eqnarray} \label{fix2}
(X^n)^g/Z_g \cong \prod_{l =1}^n X^{N_l} /S_{N_l} = \prod_{l=1}^n X^{(N_l)}.
\end{eqnarray}

\subsection{Reduction to cycles} \label{sec:reductiontocycles}
For $g \in S_n$ of cycle type $(N_1, \cdots, N_n)$,
it is easy to see that 
\begin{eqnarray*}
\chi_g(X^n, E^{\boxtimes n}) 
= \prod_{l=1}^n \chi_{\sigma_l}(X^l, E^{\boxtimes l})^{N_l},
\end{eqnarray*}
where $\sigma_l = (12\cdots l)$.
Therefore,
we have
\begin{eqnarray*}
&& \sum_{n \geq 0} \chi(X^n/S_n, E^{\boxtimes n}/S_n) p^n 
= \sum_{n \geq 0} p^n \sum_{\sum lN_l = n} 
\frac{1}{\prod_{l=1}^n N_l! l^{N_l}}
\prod_{l=1}^n (\chi_{\sigma_l} (X^l, E^{\boxtimes l}))^{N_l} \\
& = & \sum_{n \geq 0} \sum_{\sum lN_l = n} \prod_{l=1}^n
\frac{1}{N_l! l^{N_l}}
(p^l \chi_{\sigma_l} (X^l, E^{\boxtimes l}))^{N_l}
= \prod_{l \geq 1} \sum_{N_l \geq 0} \frac{1}{N_l! l^{N_l}}
(p^l \chi_{\sigma_l} (X^l, E^{\boxtimes l}))^{N_l} \\
& = & \prod_{l \geq 1} \exp 
\left(\frac{1}{l} p^l \chi_{\sigma_l} (X^l, E^{\boxtimes l}) \right). 
\end{eqnarray*}
To summarize, we have 
\begin{eqnarray}\label{eqn:reductiontocycles1}
\sum_{n \geq 0} \chi(X^n/S_n, E^{\boxtimes n}/S_n) p^n 
= \prod_{l \geq 1} \exp 
\left(\frac{1}{l} p^l \chi_{\sigma_l} (X^l, E^{\boxtimes l}) \right).
\end{eqnarray}

\subsection{Calculations for cycles}
\label{sec:calculationsforcycles}

First notice that $TX^n|_{\Delta_n(X)}$ 
is isomorphic to the direct sum of $n$ copies of $TX$ and
$\sigma_n$ acts on it by cycling the factors.
Since the eigenvalues of an $n$-cycle on $\bC^n$ are $e^{2 k \pi \sqrt{-1}/n}$
for $k = 0, \cdots, n-1$, 
each with multiplicity $1$,
it is clear that
$$\phi_{\sigma_n}(TX^n|_{\Delta_n(X)}) \cong
\sum_{k = 0}^{n-1} e^{2 k \pi \sqrt{-1}/n} TX,$$
where $k = 0$ corresponds to $T\Delta_n(X)$.
Hence 
$$ \phi_{\sigma_n}(N_{\Delta_n(X)/X^n}) \cong
\sum_{k = 1}^{n-1} e^{2 k \pi \sqrt{-1}/n} TX.$$
If $\{x_j: j =1, \cdots, d\}$ are the Chern roots of $TX$
($d = \dim X$),
then we have
\begin{eqnarray*}
\ch_{\sigma_n} \Lambda_{-1} ((N_{\Delta_n(X)/X^n})^*)
= \prod_{j=1}^d \prod_{k=1}^{n-1} (1 - e^{-x_j - 2 k \pi\sqrt{-1}/n}) 
= \prod_{j=1}^d \frac{1 - e^{-nx_j}}{1 - e^{-x_j}}.
\end{eqnarray*}
This gives
\begin{eqnarray} \label{ncyclenormal}
\ch_{\sigma_n} \Lambda_{-1} ((N_{\Delta_n(X)/X^n})^*)
= \frac{n^d \cT (TX)}{\cT(\psi^n(TX))}.
\end{eqnarray}

Combining (\ref{Lefschetz2}), (\ref{ncyclenormal}) 
and Corollary \ref{cor:Adams},
we get
\begin{eqnarray} \label{chisigman}
\chi_{\sigma_n}(X^n, E^{\boxtimes n}) 
= \frac{1}{n^d} \int_X \ch(\psi^n(E)) \cT(\psi^n(TX)).
\end{eqnarray}

\begin{lemma} \label{lm:chisigman}
We have $\chi_{\sigma_n}(X^n, E^{\boxtimes n}) = \chi(X, E)$.
\end{lemma}

\begin{proof}
Denote by $\{x_1, \cdots, x_d\}$ and $\{y_1, \cdots, y_r\}$ 
the Chern roots of $TX$ and $E$ respectively,
then we have
\begin{eqnarray*}
\chi(X, E) 
= \int_X \sum_{i=1}^r e^{y_i} \cdot \prod_{j=1}^d \frac{x_j}{1 - e^{-x_j}}
= \int_X F_d(x_1, \cdots, x_d, y_1, \cdots, y_r)
\end{eqnarray*}
for some homogeneous polynomial $F_d$ of degree $d$.
Now from (\ref{chisigman})
we have
\begin{eqnarray*}
&& \chi_{\sigma_n}(X^n, E^{\boxtimes n}) 
= \frac{1}{n^d} \int_X \sum_{i=1}^r e^{n y_i} \cdot \prod_{j=1}^d 
\frac{nx_j}{1 - e^{-nx_j}} \\
& = & \frac{1}{n^d} \int_X F_d(nx_1, \cdots, nx_d, ny_1, \cdots, ny_r) \\
& = & \int_X F_d(x_1, \cdots, x_d, y_1, \cdots, y_r) 
  = \chi(X, E).
\end{eqnarray*}
\end{proof}

\begin{theorem}
For a closed complex manifold $M$ and a holomorphic vector bundle $E$ on $X$,
we have
\begin{eqnarray*}
\sum_{n \geq 0} \chi(X^n/S_n, E^{\boxtimes n}/S_n) p^n 
= \frac{1}{(1 - p)^{\chi(X, E)}}, 
\end{eqnarray*}
\end{theorem}

\begin{proof}
Combining (\ref{eqn:reductiontocycles1}) with Lemma \ref{lm:chisigman},
we get
\begin{eqnarray*}
&& \sum_{n \geq 0} \chi(X^n/S_n, E^{\boxtimes n}/S_n) p^n 
=  \prod_{l \geq 1} \exp 
\left(\frac{1}{l} p^l \chi_{\sigma_l} (X^l, E^{\boxtimes l}) \right) \\
& = & \prod_{l \geq 1} \exp 
\left(\frac{1}{l} p^l \chi(X, E)\right) = \frac{1}{(1 - p)^{\chi(X, E)}}. 
\end{eqnarray*}
\end{proof}

\begin{remark}
One can find an easy proof of the formula in the above theorem using
$H^*(X^n/S_n, E^{\boxtimes n}/S_n) \cong S^n(H^*(X, E))$.
\end{remark}

\subsection{Formula for $\chi(X^n, E^{\boxtimes n}|S_n)$}
Similarly, we have
\begin{eqnarray*}
&& \sum_{n \geq 0} \chi (X^n, E^{\boxtimes n}| S_n) p^n \\
& = & \sum_{n \geq 0} p^n \sum_{\sum lN_l = n} \prod_{l=1}^n 
\chi (X^{N_l}/S_{N_l}, 
(\phi_{\sigma_l}(E^{\boxtimes l}))^{\boxtimes N_l}/S_{N_l}) \\
& = & \sum_{n \geq 0} \sum_{\sum lN_l = n} \prod_{l=1}^n p^{lN_l}
\chi (X^{N_l}/S_{N_l}, 
(\phi_{\sigma_l}(E^{\boxtimes l}))^{\boxtimes N_l}/S_{N_l}) \\
& = & \prod_{l \geq 1} \sum_{N_l \geq 0} 
p^{N_l l} \chi (X^{N_l}/S_{N_l}, 
(\phi_{\sigma_l}(E^{\boxtimes l}))^{\boxtimes N_l}/S_{N_l}) \\
& = & \prod_{l \geq 1} \prod_{m \geq 1} 
\exp \left( \frac{1}{m} p^{lm} 
\chi_{\sigma_m}(X^m, (\phi_{\sigma_l}(E^{\boxtimes l}))^{\boxtimes m}) \right).
\end{eqnarray*}
To summarize,
we have
\begin{eqnarray} \label{eqn:reductiontocycles2} 
\;\;\;\;\;\;
\sum_{n \geq 0} \chi (X^n, E^{\boxtimes n}| S_n) p^n 
= \prod_{l \geq 1} \prod_{m \geq 1} 
\exp \left( \frac{1}{m} p^{lm} 
\chi_{\sigma_m}(X^m, (\phi_{\sigma_l}(E^{\boxtimes l}))^{\boxtimes m}) \right).
\end{eqnarray}

\begin{theorem}
For a closed complex manifold $M$ and a holomorphic vector bundle $E$ on $X$,
we have
\begin{eqnarray*}
\sum_{n \geq 0} \chi(X^n, E^{\boxtimes n}|S_n) p^n 
= \prod_{l \geq 1}\frac{1}{(1 - p^l)^{\chi(X, \psi^k(E))}}. 
\end{eqnarray*}
\end{theorem}

\begin{proof}
Combining (\ref{eqn:reductiontocycles2}) with Proposition \ref{prop:Adams} 
and Lemma \ref{lm:chisigman} for $\psi^{l}(E)$,
we get
\begin{eqnarray*}
&& \sum_{n \geq 0} \chi (X^n, E^{\boxtimes n}| S_n) p^n 
= \prod_{l \geq 1} \prod_{m \geq 1} 
\exp \left( \frac{1}{m} p^{lm} 
\chi_{\sigma_m}(X^m, (\phi_{\sigma_l}(E^{\boxtimes l}))^{\boxtimes m}) \right) \\
& = & \prod_{l \geq 1} \prod_{m \geq 1} \exp \left( \frac{1}{m} p^{lm} 
\chi_{\sigma_m}(X^m, (\psi^l(E))^{\boxtimes m}) \right) \\
& = & \prod_{l \geq 1} \prod_{m \geq 1} \exp \left( \frac{1}{m} p^{lm} 
\chi (X, \psi^l(E)) \right) 
= \prod_{l \geq 1} \exp \left( \sum_{m \geq 1} \frac{1}{m} p^{lm} 
\chi (X, \psi^l(E)) \right) \\
& = &  \prod_{l \geq 1}\frac{1}{(1 - p^l)^{\chi(X, \psi^k(E))}}.
\end{eqnarray*}
\end{proof}
\section{Generalizations to the graded vector bundles}

\subsection{Graded anti-symmetric $S_n$-action}

We now assume that $\pi: E \to X$ is $\bZ$ or $\bZ_2$-graded vector bundle.
Define the graded anti-symmetric 
$S_n$-action on $E^{\boxtimes n}$ as follows:
\begin{eqnarray*}
\sigma^a(v_1 \boxtimes \cdots \boxtimes v_n)
= (-1)^{\epsilon(\sigma, |v_1|, \cdots, |v_n|)} 
v_{\sigma(1)} \boxtimes \cdots \boxtimes v_{\sigma(n)},
\end{eqnarray*}
where the sign $(-1)^{\epsilon(\sigma, |v_1|, \cdots, |v_n|)}$ is
determined as follows: use transpositions of adjacent vectors to change
$v_{\sigma(1)} \boxtimes \cdots \boxtimes v_{\sigma(n)}$ back to 
$v_1 \boxtimes \cdots \boxtimes v_n$,
for each such transposition $v \otimes w \mapsto w \otimes v$,
introduce a sign $(-1)^{|v||w|}$.
One can easily check that the final result does not depends 
on the choices of the transpositions.
Simialarly define the graded anti-symmetric $S_n$-action on $E^{\otimes n}$.
These definitions are motivated by the following

\begin{example} \label{exm:graded}
Let $E = \Lambda^*(T^*X)$.
Each fiber of $\Lambda^*(T^*X^n)$ is spanned by elements of the form
$p_1^*\alpha_1 \wedge \cdots \wedge p_n^*\alpha_n$,
where as earlier $p_j: X^n \to X$ is the projection onto the $j$-th factor,
$\alpha_1, \cdots, \alpha_n \in \Lambda^*(T^*X)$.
Therefore we have an isomorphism of vector bundles
\begin{eqnarray} \label{exterior}
\Lambda^*(T^*X^n) \cong \Lambda^*(T^*X)^{\boxtimes n}, \;\;\;
p_1^*\alpha_1 \wedge \cdots \wedge p_n^*\alpha_n \mapsto
\alpha_1 \boxtimes \cdots \boxtimes \alpha_n.
\end{eqnarray}
The $S_n$-action on $X^n$ induces an $S_n$-action on $\Lambda^*(T^*X^n)$:
\begin{eqnarray*}
&& \sigma^*(p_1^*\alpha_1 \wedge \cdots \wedge p_n^*\alpha_n)
= p_{\sigma(1)}^*\alpha_1 \wedge \cdots \wedge p_{\sigma(n)}^*\alpha_n \\
& = & (-1)^{\epsilon(\sigma, |\alpha_1|, \cdots, |\alpha_n|}
p_1^*\alpha_{\sigma^{-1}(1)} \otimes \cdots \otimes 
p_n^*\alpha_{\sigma^{-1}(n)}.
\end{eqnarray*}
With respect to the isomorphism (\ref{exterior}),
the induced $S_n$-action on $\lambda^*(T^*X^n)$ is just the
graded anti-symmetric action of $S_n$ on $\Lambda^*(T^*X)$.
Similar discussions can be carried out for $E = \Lambda^*(TX)$.
\end{example}

\subsection{Graded $K$-theory and graded Adams operations}

It is straightforward to define the graded $K$-theory $GK(X)$
as the Grothendieck algebra of graded vector bundles on $X$.
Similarly define the graded equivariant $K$-theory.
Define the graded Adams operation $G\psi^n: GK(X) \to GK(X)$ as follows:
if $E = L_1 \oplus \cdots \oplus L_r$,
where $L_1, \cdots, L_r$ are line bundles of degree $d_1, \cdots, d_r$ 
respectively,
then
$$G\psi^n(E) = \sum_{j=1}^r (-1)^{(n-1)d_j} L_j^{\otimes n}.$$
As in the ordinary case,
we have maps
\begin{eqnarray*}
\otimes n: GK(X) \to GK_{S_n}(X), 
\boxtimes n: GK(X) \to GK_{S_n}(X^n)
\end{eqnarray*}
defined by sending $E$ to $E^{\otimes n}$ and $E^{\boxtimes n}$ with the graded
anti-symmetric $S_n$-actions respectively.
The graded version of Proposition \ref{prop:Adams} and Corollary \ref{cor:Adams}
is the following

\begin{proposition} \label{prop:gradedAdams}
For a graded vector bundle $E$ 0n $X$,
we have
\begin{eqnarray*}
&& \phi_{\sigma_n} (E^{\otimes n}) = G\psi^n(E), \\
&& i^*\phi_{\sigma_n} (E^{\boxtimes n} = G\psi^n(E).
\end{eqnarray*}
\end{proposition}

\begin{proof}
By splitting principle,
we may assume that $E = \oplus_{m = 1}^r L_m$,
where each $L_m$'s is a graded line bundle of degree $d_m$.
Then
$$E^{\otimes n} = \bigoplus_{1 \leq j_1, \cdots, j_n \leq r}  
L_{j_1} \otimes \cdots \otimes L_{j_n}.$$
For each multiple index $J = (j_1, \cdots, j_n)$,
agian let $L_J = L_{j_1} \otimes \cdots \otimes L_{j_n}$,
$V_J$ the subbundle spanned by
line bundles $L_{\sigma_n^k(J)}$, $k=0, \cdots, n-1$.
Then $V_J$ is invariant under the action of $\bZ_n = \langle \sigma_n \rangle$.
Since $E^{\otimes n}$ is a direct sum of such $V_J$,s, 
it suffices to find $\phi_{\sigma_n}(V_J)$ for all $J$.
When $j_1 = \cdots = j_n = m$,
$V_J = L_m^{\otimes n}$ and $\sigma_n$ acts by the multiplication of
$$(-1)^{\epsilon(\sigma_n, d_m, \cdots, d_m)} 
= (-1)^{d_m \cdot (n-1)d_m} = (-1)^{(n-1)d_m}.$$
In other words,
$$\phi_{\sigma_n}(\sum_{m=1}^r y^{nd_m} L_m^{\otimes n}) 
= \sum_{m=1}^r (-1)^{(n-1)d_m} (y^{d_m} L_m)^{\otimes n} 
= G\psi^n(E).$$
When $j_1, \cdots, j_n$ are not all identical,
we modify the proof of Proposition \ref{prop:Adams} by computing 
the characters as explained in Remark \ref{remark:character} to show 
$\phi_{\sigma_n}(V_j) = 0$.
\end{proof}

\subsection{$\chi_{-y}$ of graded holomorphic vector bundles}

For a holomorphic $\bZ$-graded vector bundle $E = \oplus_{j=0}^m E^j$ on $X$,
write $E_{-y} = \sum_{j=0}^m (-y)^jE^j$ and
$$\chi_{-y}(X, E) = \chi(X, E_{-y}) = \sum_{j=0}^m (-y)^j \chi(X, E).$$
There is an induced grading on $E^{\boxtimes E}$ and we have
$$E^{\boxtimes n}_{-y} = (E_{-y})^{\boxtimes n}$$
in the sense that both sides can be written as
$$\sum_{k \geq 0} (-y)^k \sum_{j_1  + \cdots j_n = k} 
E^{j_1} \boxtimes \cdots \boxtimes E^{j_n}.$$
Denote by $E^{\boxtimes n}/S_n^a$ the quotient of $E^{\boxtimes n}$
by the graded anti-symmetric action.
Now we have
\begin{eqnarray*}
\chi_{-y}(X^n/S_n, E^{\boxtimes n}/S_n^a)
= \sum_{[g] \in (S_n)_*} \frac{1}{|Z_g|} \chi_{g^a}(X^n, E^{\boxtimes n}_y).
\end{eqnarray*}
Here we have used $g^a$ to indicated $g$ acts graded anti-symmetrically.

\begin{theorem} \label{thm:graded}
Let $E$ be a $\bZ$-graded holomorphic vector bundle over a closed 
complex manifold $X$,
then we have
\begin{eqnarray*}
\sum_{n \geq 0} \chi_{-y}(X^n/S_n, E^{\boxtimes n}/S^a_n) p^n 
= \prod_{l \geq 1} \exp 
\left(\frac{1}{l} p^l \chi_{-y^l} (X, E) \right).
\end{eqnarray*}
\end{theorem}

\begin{proof}
As in \S \ref{sec:reductiontocycles},
we have
\begin{eqnarray*}
\sum_{n \geq 0} \chi_{-y}(X^n/S_n, E^{\boxtimes n}/S^a_n) p^n 
= \prod_{l \geq 1} \exp 
\left(\frac{1}{l} p^l \chi_{\sigma_l^a} (X^l, E_{-y}^{\boxtimes l}) \right).
\end{eqnarray*}
From (\ref{Lefschetz2}), (\ref{ncyclenormal}) 
and Proposition \ref{prop:gradedAdams},
we get
\begin{eqnarray*} 
\chi_{\sigma_n^a}(X^n, E_{-y}^{\boxtimes n}) 
= \frac{1}{n^d} \int_X \ch(G\psi^n(E_{-y})) \cT(\psi^n(TX)).
\end{eqnarray*}
Without loss of generality,
we may assume that $E = \oplus_{m = 1}^r L_m$,
where each $L_m$'s is a graded line bundle of degree $d_m$.
Let $y_m = c_1(L_m)$.
Also let $\{x_1, \cdots, x_d\}$ be the Chern roots of $TX$.
Since
$$G\psi^n(E_{-y}) = G\psi^n(\sum_{i=1}^r (-y)^{d_i} L_i) = 
= \sum_{i=0}^r (-1)^{(n-1)d_i} (-y)^{nd_i} L_i^{\otimes n},$$
we have
\begin{eqnarray*} 
\chi_{\sigma_n^a}(X^n, E_{-y}^{\boxtimes n}) 
& = & \frac{1}{n^d} \int_X \sum_{i=0}^r (-1)^{(n-1)d_i}(-y)^{nd_i}e^{ny_i}
\prod_{j=1}^d \frac{x_j}{1 - e^{-nx_j}} \\
& = & \int_X \sum_{i=0}^r (-1)^{(n-1)d_i}(-y)^{nd_i}e^{y_i}
\prod_{j=1}^d \frac{x_j}{1 - e^{-x_j}} \\
& = & \int_X \ch(E_{-y^n})) \cT(TX) =\chi(X, E_{-y^n}) = \chi_{-y^n}(X, E).
\end{eqnarray*}
Here for the second equality we have used the same argument as in the proof of
Lemma \ref{lm:chisigman}.
This completes the proof.
\end{proof}

\subsection{The Hirzebruch $\chi_{-y}$-genera of the symmetric products}

Recall Hirzebruch's $\chi_{-y}$-genus for a compact complex manifold $M$ is
$$\chi_{-y}(M) = \sum_{p, q \geq 0} (-1)^qh^{p, q}(M) (-y)^p 
= \chi_{-y}(M, \Lambda^*(T^*M)) =\chi(M, \Lambda_{-y}(T^*M)).$$
By Example \ref{exm:graded}, 
\begin{eqnarray*} 
\sum_{n \geq 0} \chi_{-y}(X^{(n)}) q^n 
= \sum_{n \geq 0} \chi(X^n/S_n, \Lambda_{-y}(T^*X^n)/S^a_n) q^n.
\end{eqnarray*}
This recovers (\ref{gfchiy1}) by Theorem \ref{thm:graded}.

\subsection{The $B$-Hirzebruch genera $\hat{\chi}_{-y}$ of symmetric products}
As explained in Zhou \cite{Zho},
the study of mirror symmetry of Calabi-Yau manifolds
motivates the following definition:
$$\hat{\chi}_{-y}(M) = \chi_{-y}(M, \Lambda^*(TM)) 
	= \chi(M, \Lambda_{-y}(TM)).$$
Denote by $h^{-p, q}(M)$ the dimension of $H^q(M, \Lambda^p(TM))$,
then
$$\hat{\chi}_{-y}(M) = \sum_{p, q \geq 0} (-1)^qh^{-p, q}(M)(-y)^p.$$
We have proved in \cite{Zho} the following formula:
\begin{eqnarray} \label{gfhatchiy} 
\sum_{n \geq 0} \hat{\chi}_{-y}(X^{(n)}) q^n 
= \exp \left(\sum_{l \geq 1} \frac{\hat{\chi}_{-y^l}(X)}{l} q^l\right).
\end{eqnarray}
This can be recovered by exactly the same method 
(Example \ref{exm:graded} and Theorem \ref{thm:graded})
as for $\chi_{-y}$.

\subsection{Formulas for $\chi_{-y}(X^n, S_n)$ and $\hat{\chi}_{-y}(X^n, S_N)$}
From (\ref{complexreducetofix2}) and (\ref{fix2}),
one gets:
\begin{eqnarray*}
&& \sum_{n \geq 0} \chi_{-y}(X^n, S_n) q^n
= \sum_{n \geq 0} q^n \sum_{\sum lN_l = n} 
\prod_{l=1}^n y^{N_l(l-1)d/2}\chi_{-y}(X^{(N_l)}) \\
& = & \prod_{l \geq 1} \sum_{N_l \geq 0} 
	q^{lN_l}y^{N_l(l-1)d/2}\chi_{-y}(X^{(N_l)}) 
= \prod_{l \geq 1} \sum_{N_l \geq 0} 
	(q^ly^{(l-1)d/2})^{N_l}\chi_{-y}(X^{(N_l)})\\
& = & \prod_{l \geq 1} \exp \left( 
\sum_{n \geq 1} \frac{\chi_{-y^n}(X)}{n} (q^ly^{(l-1)d/2})^n\right) 
\;\;\;\;\;\;\;\; \mbox{by (\ref{gfchiy1})} \\
& = & \exp \left( \sum_{n \geq 1} \frac{p^n \chi_{-y^n}(X)}{n}
	\sum_{l \geq 1} (q^ny^{nd/2})^{l-1}\right)
= \exp \left( \sum_{n \geq 1} \frac{q^n}{n} 
\frac{\chi_{-y^n}(X)}{1 - (y^{d/2}q)^n} \right).
\end{eqnarray*}
This recovers (\ref{gfchiy2}).
Similarly,
we have
$$\hat{\chi}_y(M) = \sum_{[g] \in G_*} (-y)^{F_g} \hat{\chi_y}(M^g/Z_g),$$
for a complex orbifold $M/G$.
Then the same proof recovers:
\begin{eqnarray*}
\sum_{n \geq 0} \hat{\chi}_{-y}(X^n, S_n) q^n
= \exp \left( \sum_{n \geq 1} \frac{q^n}{n} 
\frac{\hat{\chi}_{-y^n}(X)}{1 - (y^{d/2}q)^n} \right).
\end{eqnarray*}

{\bf Acknowledgement}. 
{\em The work in this paper began when the author was visiting
the Morningside Mathematics Center in Beijing in the summer of 1999.
It was finished during the author's visit at Department of Mathematics,
Texas A\&M University.
The author thanks Kefeng Liu, Xiaonan Ma, Weiqiang Wang and Weiping Zhang 
for helpful discussions.}

\end{document}